\documentclass[12pt,reqno]{article}

\usepackage{amssymb,amsfonts,amsmath}
 
\usepackage[a4paper]{geometry}

\usepackage{amssymb}
 
\usepackage{parskip}
 
\usepackage{graphicx}

\usepackage{textcomp}
 
\newtheorem{theorem}{Theorem}[section]
\newtheorem{lemma}[theorem]{Lemma}

\newtheorem{corollary}[theorem]{Corollary}

 \newtheorem{note}[theorem]{Note}
\usepackage{mathtools}
\DeclarePairedDelimiterX{\norm}[1]{\lVert}{\rVert}{#1}

\usepackage{fancyhdr}
\title{Rotation Sets of Billiards with N Obstacles on a Torus }
\author{Zainab Alsheekhhussain}
\begin{document}
\def\endofproof{{\rule{6pt}{6pt}}}
\def\i{{\bf{i}}}
\def\j{{\bf{j}}}
\def\0{{\bf{0}}}
\def\l{{\bf{l}}}
\def\x{{\bf{x}}}
\def\y{{\bf{y}}}
\def\u{{\bf{u}}}
\def\v{{\bf{v}}}
\def\n{{\bf{n}}}
\def\k{{\bf{k}}}

\titlepage
\maketitle

\begin{abstract}
For billiards with $N$ obstacles on a torus, we study the behavior of specific kind of its trajectories, \emph{the so called admissible trajectories}. Using the methods developed in \cite{1}, we prove that the \emph{admissible rotation set} is convex, and the periodic trajectories of admissible type are dense in the admissible rotation set. In addition, we show that the admissible rotation set is a proper subset of the general rotation set. 

\end{abstract}


\let\thefootnote\relax\footnotetext{School of Mathematics and Statistics\\ University of Western Australia, Perth WA 6009\\ Australia\\zainab.alsheekhhussain@research.uwa.edu.au}

\section{Introduction}
Rotation Theory plays an important role in the study of various dynamical systems. \emph{Rotation vectors} and \emph{sets} provide essential information for understanding the behavior of trajectories of a billiard system \cite{1}, \cite{2}.

The dynamical system we investigate in this paper is the billiard on an $m$-dimensional torus in the exterior of $N$ convex obstacles. We associate with it a natural observable--the so called \emph{displacement} as we consider continuous time. The main considerations are on the universal covering $\mathbb{R}^m$ of the torus.

 The methodology that we use in this paper has been developed by Blokh, Misiurewicz and Sim\'{a}nyi in \cite{1}, where they have studied billiard on a torus with one convex obstacle. In the present paper, we extend their arguments and show how to deal with an arbitrary number of convex obstacles satisfying a certain no-eclipse condition.

This model is very important in statistical physics as it is related to the hard ball model and the Lorentz gas model (see \cite{3}, \cite{4}).

To summarize our results, we have to emphasize first that we deal mainly with a specific kind of trajectories--the so called \emph{admissible trajectories}--that satisfy certain properties. We prove that the set consisting of rotation vectors of those kind of trajectories is convex. In addition, we prove that the rotation vectors of periodic trajectories of admissible type are dense in this set. 

However, whether the general rotation set is convex, or the periodic trajectories are dense in it are two natural questions that need a separate investigation. We believe that the results in this paper could help in this direction. 

We refer to reader to \cite{4}, \cite{13}, \cite{7} and \cite{8} for the definition and general information about billiards and billiards on a torus in particular.

Concerning general information about Rotation Theory and references, the reader could see \cite{13}, \cite{12}, \cite{28}, \cite{24}, \cite{25} and  \cite{11}.

\section{Billiards on a torus}
\renewcommand{\theequation}{\arabic{section}.\arabic{equation}}
Consider the class of billiards on the $m$-dimensional torus $\mathbb{T}^m=\mathbb{R}^m/\mathbb{Z}^m$ $(m\geq 2)$ in the exterior of $N$ obstacles $O_r$, $r=1,...,N$, where the obstacles are strictly convex disjoint bodies with smooth $(C^2)$ boundaries  of diameter less than $ \frac{\sqrt 2}{4}$.

Let $\phi:\mathbb{R}^m\longrightarrow\mathbb{T}^m$ be the natural covering of $\mathbb{T}^m$ where $\phi([0,1]^m)=\mathbb{T}^m$, and let $$\phi^{-1}(O_r)\cap[0,1]^m=O_{(\0, r)}.$$
That is to say that we assume the covering and the coordinates in $\mathbb{T}^m$ are chosen so that $O_{(\0, r)}$ is in the interior of $[0,1]^m$ for all $r=1,...,N$. This assumption follows from the first main assumption bellow. Setting $O_{(\k, r)}=O_{(\0, r)}+\k$ for $\k\in\mathbb{Z}^m, \, r=1,...,N$, we have $\phi^{-1}(O_r)=\cup_{\k\in\mathbb{Z}^m}O_{(\k, r)}$. We will consider billiard trajectories in $\Omega=\overline{\mathbb{R}^m\backslash K}$, where $K=\cup_{\k\in\mathbb{Z}^m}\cup_{r=1}^NO_{(\k, r)}$. Clearly, the projection via $\phi$ of such trajectories are billiard trajectories in $\overline{\mathbb{T}^m\backslash\cup_{r=1}^NO_r}$.

In this paper, we make the following 
\medskip

{\bf Assumptions:}

\begin{itemize}
\item The obstacles are located in positions such that, in $S_0 = [0,1]^m$, if $\pi_i :\mathbb{R}^m\rightarrow\mathbb{R}$ are projections defined by $ \pi_i(x_1,...,x_i,...,x_m)=x_i$ then, for every $i$ and for every $r\neq r'$, we have
\begin{equation}\pi_i(O_r)\cap\pi_i(O_{r'})=\emptyset ; \end{equation}
\item  The obstacles in $S_0$ satisfy the standard {\it no-eclipse condition}: the convex hull of any two obstacles has no common points with any of the other obstacles.  Note that we need the no-eclipse condition to be satisfied just in $S_0$ and not in the whole plane.
\end{itemize}

\medskip
We will use some terminology from \cite{1}. For $\k,\i,\j$ and $\i\neq \j$, if the obstacle $O_{(\k, r)}$ has a common point with the convex hull of $O_{(\i, s)}\cup O_{(\j, t)}$ then we will say that $O_{(\k, r)}$ is \emph{between} them. We will say that a trajectory $T$ is of \emph{type} $(\k_n,r_n)_{n=0}^\infty$ if a continuous lifting of $T$ to $\mathbb{R}^m$ starts at $O_{(\k_0,r_0)}$, where $\k_0=0$, and its consecutive reflections are from $O_{(\k_n,r_n)}$, $n\in \mathbb{N}$. In some situations, it is more convenient to study the type of a piece of a trajectory, so then we will consider a finite sequence.  We study mostly one-sided trajectories of special kind defined below. 

Note that due to Lemma 2.1 in \cite{1}, if the trajectory is not tangent at its initial point then it has infinitely many reflections.

A sequence $\{(\k_n, r_n)\}_{n=0}^p$, where $p$ could be $\infty$, satisfying the following conditions is called \emph{admissible}:

(i) $\k_0=\0$;

(ii) $(\k_{n+1},r_{n+1})\neq(\k_n,r_n)$ for every $n$;

(iii) For every $n<p$, there is no obstacle between $O_{(\k_n,r_n)}$ and $O_{(\k_{n+1},r_{n+1})}$;

(iv) For every $n<p-1$, the obstacle $O_{(\k_{n+1},r_{n+1})}$ is not between $O_{(\k_n,r_n)}$ and $O_{(\k_{n+2},r_{n+2})}$.

We will denote the length of a trajectory piece $T$ by $\mid T\mid$ and its displacement by $d(T)$, that is the length of the difference between its initial point and end point.

We define the \emph{general rotation set } $R$ to be the set consisting of all limit points of sequences of type $(\frac{d(T_n)}{\mid T_n\mid})_{n=1}^\infty=(\frac{\y_n-\x_n}{t_n})_{n=1}^\infty$ where there is a trajectory piece $T_n$ in $\mathbb{R}^m$ from $\x_n$ to $\y_n$ of length $t_n$, and $t_n$ goes to infinity. The \emph{admissible rotation set} $AR$ is the set consisting of all limit points of trajectories of admissible type. By definition, both sets are closed and contained in the closed unit ball in $\mathbb{R}^m$, and the admissible rotation set is contained in the general rotation set.

For a full trajectory $L$ starting at $\x\in \mathbb{R}^m$, if the limit of the sequence $(\frac{d(T_n)}{\mid T_n\mid})_{n=1}^\infty$, where $d(T_n)$ the length of the difference between the starting point and the $n^{th}$ reflection point,  exists then we will call it the \emph{rotation vector} of $\x$. This limit depends only on the trajectory. (I.e. if we choose another point on that trajectory to be the initial point, the limit will not change.) Hence, we can call it the \emph{rotation vector of the trajectory}.

For a periodic orbit, denote by $P$  the piece of the periodic trajectory that represents its period. Clearly, the rotation vector for any periodic trajectory exists and equals $\frac{d(P)}{\mid P\mid}$.

Now we can examine the connection between  admissible sequences and trajectories.

As in \cite{1}, we have the following.
\begin{theorem}
\label{1}
For every admissible sequence $(\k_n,r_n)_{n=0}^{\infty}$, there is a billiard trajectory of type $(\k_n,r_n)_{n=0}^{\infty}$.
\end{theorem}
\emph{Proof}. Fix $n$ and consider the map $F:Z_0\times Z_1\times ....\times Z_n\longrightarrow \mathbb{R}$, where $Z_j$ is the set consisting of all points on the boundary of the obstacle $O_{(\k_j, r_{j})}$, defined by
$$ F(\x_0^{(n)},\x_1^{(n)},.....,\x_n^{(n)})=\parallel \x_0^{(n)}-\x_1^{(n)}\parallel+\parallel \x_1^{(n)}-\x_2^{(n)}\parallel+...+\parallel \x_{n-1}^{(n)}-\x_n^{(n)}\parallel,$$ $\x_j^{(n)}$ belongs to the boundary of the obstacle $O_{(\k_j, r_{j})}$. Clearly, $F$ is continuous and the value of $F$ depends only on the choice of $\x_j^{(n)}$. In addition, the space $Z_0\times Z_1\times ...\times Z_n$ is compact. Hence, $F$ has an absolute minimum at $(\x_0^{(n)},\x_1^{(n)},...,\x_n^{(n)})$.

Let $\eta$ be the minimum curve joining $\x_j$'s by straight segments. By property (iii) of an admissible sequence, we have that any segment joining two consecutive points $\x_j^{(n)}, \x_{j+1}^{(n)}$ cannot intersect any obstacle other than $O_{(\k_j, r_{j})}, O_{(\k_{j+1}, r_{j+1})}$ (the ones they belong to). By minimality of $\eta$ and property (iv), we have that the segment cannot intersect $O_{(\k_j, r_{j})}, O_{(\k_{j+1}, r_{j+1})}$ at more than one point. Now $\eta$ is minimal, so at every $\x_j^{(n)}$, $j=1,...,n-1$, the angle of incidence is equal to the angle of reflection. Thus, $\eta$ is a piece of trajectory.

By the argument above, for $j$ greater or equal to $1$, for any admissible sequence $(\k_n,r_n)_{n=0}^j$, we can find a trajectory piece $\eta=(\x_0^j, \x_1^j,....,\x_j^j)$ of that type. Let the following be trajectory pieces:

\begin{eqnarray}
&& \x_0^{(1)}, \x_1^{(1)}\nonumber\\
&&\x_0^{(2)}, \x_1^{(2)}, \x_2^{(2)}\nonumber\\
&&\x_0^{(3)}, \x_1^{(3)}, \x_2^{(3)}, \x_3^{(3)}\\
&&..................................\nonumber\\
&&\x_0^{(n)}, \x_1^{(n)},........., \x_n^{(n)}\nonumber\\
&&..................................\nonumber
\end{eqnarray}

Looking at the first vertical sequence $$\x_0^{(1)},\x_0^{(2)},.....,\x_0^{(n)},.... ,$$ we can choose a convergent subsequence, e.g. $$\x_0^{(k_1^{(0)})}, \x_0^{(k_2^{(0)})},.....,\x_0^{(k_n^{(0)})},....$$ converges to some $\x_0$. Then look at the corresponding subsequence in the second column of (2.2). This is $$\x_1^{(k_1^{(0)})},\x_1^{(k_2^{(0)})},.....,\x_1^{(k_n^{(0)})},..... .$$ Choose a convergent subsequence of this, e.g. $$\x_1^{(k_1^{(1)})},\x_1^{(k_2^{(1)})},.....,\x_1^{(k_n^{(1)})},.....$$ converges to some $\x_1$. Continuing by induction at the $p^{th}$ step, we construct a convergent subsequence $$\x_p^{(k_1^{(p)})}, \x_p^{(k_2^{(p)})},.....,\x_p^{(k_n^{(p)})},....$$ of the $p^{th}$ column of (2.2) converging to some $\x_p$ . Then the sequence $(\x_i)_{i=0}^\infty$ is the required trajectory. \endofproof

Note that the above Theorem can be generalized to the case where the admissible sequence is of the form $(\k_n,r_n)_{-\infty}^{\infty}$. In addition, for an admissible sequence where $(\k_{n+q}, r_n)=(\k_{n+p}, r_n)$ for every $n$, where $p$ is in $\mathbb{Z}^m$ and $q$ is a positive integer, we can choose a periodic trajectory of this type with discrete period $q$. 

\begin{note}\label{2} Lemma 2.3 and Corollary 2.4 from \cite{1} show that if $(\k_n,r_n)_{n=0}^s$ is a finite sequence  (not necessary admissible) of elements of $\mathbb{Z}^m\times \{1,...,N\}$, and if $\x_0$ is in $\partial O_{(\k_0, r_0)}$ and $\x_s$ is in $\partial O_{(\k_s, r_s)}$ then there is at most one trajectory piece of this type starting at $\x_0$ and ending at $\x_s$ even when the first segment of the trajectory piece intersects $O_{(\k_0, r_0)}$ and the last segment intersects $O_{(\k_n, r_n)}$. Moreover, if the trajectory piece has admissible type then it is the shortest path of that type starting and ending at those points.
\end{note}

Let $d_r$ be {\it{the diameter of the obstacle}} $O_r$ for every $r=1,...,N$, and let $d$ be {\it{the maximal}} $d_r$. The next lemma will show that, for any finite admissible sequence, the lengths of the trajectory pieces of its type depend on the points $\x_0$ and $\x_s$, but they can differ by at most $2d$.

\begin{lemma}
\label{3}
The difference between the lengths of the trajectory pieces of the same admissible sequence of type $(\k_n,r_n)_{n=0}^s$ is at most $2d$.
\end{lemma}
\emph{Proof}. Following the procedure in \cite{1}, take any two trajectory pieces of type $(\k_n,r_n)_{n=0}^s$ starting at $\x_0$, $\y_0$ and ending at $\x_s$, $\y_s$ respectively. Then, by adding the segments connecting $\x_0$ with $\y_0$ and $\x_s$ with $\y_s$ at the beginning and end of the first trajectory piece, we get a path joining $\y_0$ with $\y_s$ of length not shorter than the length of the second trajectory piece, by Note \ref{2}, but can exceed its length by at most $d_{r_0}+d_{r_s}\leq 2d$. As the two trajectory pieces are arbitrary, we prove the lemma. \endofproof
\begin{note}
\label{4}
Clearly, the difference between the displacements of trajectory pieces of the same admissible sequence type is at most $2d$.
\end{note}
\section{Directed Graph}
\setcounter{equation}{0}
Again we will use an idea from \cite{1} with appropriate modifications.

In order to define our directed graph $G$, we need to look at the definition of an admissible sequence from different prospective. For a sequence  $(\k_n, r_n)$ of ordered pairs of elements of $\mathbb{Z}^m\times \{1,....,N\}$, let $(\l_n, r_n)=(\k_n-\k_{n-1}, r_n)$ and $(\l_0, r_0)=(\0, r_0)$. Then the sequences $ (\k_n, r_n)_{n=1}^\infty$ and $ (\l_n, r_n)_{n=1}^\infty$ are equivalent, i.e. from each of them we can recover the other.  Hence, we can restate conditions (iii) and (iv) as follows:
(iii) For every $n$, there is no obstacle between $O_{(\0,r_{n-1})}$ and $O_{(\l_n, r_n)}$;
(iv) The obstacle $O_{(\l_n, r_n)}$ is not between $O_{(\0, r_{n-1})}$ and $O_{(\l_n+\l_{n+1}, r_{n+1})}$.

\medskip

{\bf{Definition of Directed Graph:}}

The vertices of the directed graph are the ordered pairs $(\j, s)$, where $\j\in \mathbb{Z}^m$, $s=1,2,....,N$ and there is no $r=1,2,...,N$ such that there is an obstacle between $O_{(\0, r)}$ and $O_{(\j, s)}$. There is an edge from $(\j, s)$ to $(\i, t)$, where $\j,\i \in \mathbb{Z}^m$ if and only if there exists $r=1,2,..,N$ such that $O_{(\j, s)}$ is not between $O_{(\0, r)}$ and $O_{(\i+\j, t)}$.

\medskip
Hence, we can deal with admissible sequences $(\l_n, r_n)_{n=1}^\infty$ as one-sided paths in $G$. Note that there is no vertex $(\j, s)$ with an edge from it to itself. In addition, $G$ is symmetric. That is, if $(\j,s)$ is a vertex then $(-\j,s)$ is a vertex; there is an edge from $(\j,s)$ to $(-\j,s)$, and if there is an edge from $(\j, s)$ to $(\k, t)$ then there is an edge from $(-\j, s)$ to $(-\k, t)$. 

Note that, by our assumption (the obstacles in $S_0$ satisfy the no-eclipse condition), all pairs $(\0,r)$, $r=1,2,...,N$, are vertices and there is an edge between any two of them.

Lemma 2.7 in \cite{1} shows another kind of symmetry in $G$ in the case of one obstacle, and we restate this lemma with some modifications to be able to apply it to our case.
\begin{lemma}
\label{5}
Let $\k, \j$ be in $\mathbb{Z}^m$. If $O_{(\k,t)}$ is not between $O_{(\0, r)}$ and $O_{(\k+\j, s)}$, then $O_{(\j,t)}$ is not between $O_{(\0, s)}$ and $O_{(\k+\j, r)}$.
\end{lemma}
\emph{Proof}. Consider the map $f(\x)=\k+\j-\x$. The image of $O_{(\0,r)}$ is $O_{(\k+\j, r)}$ and $f(O_{(\k+\j, s)})=O_{(\0, s)}$, $f(O_{(\k, t)})= O_{(\j, t)}$, which proves the lemma. \endofproof

\begin{corollary}
\label{6}
If there is an edge in $G$ from $(\k, t)$ to $(\j, s)$, so there is $r$ such that $O_{(\k, t)}$ is not between $O_{(\0, r)}$ and $O_{(\k+\j, t)}$, then there is an edge from  $(\j, t)$ to $(\k, r)$.
\end{corollary}

 We want now to show that our graph $G$ is connected, but first we will need several lemmas.
\begin{lemma}
\label{7}
The set of vertices of $G$ is finite.
\end{lemma}
\emph{Proof}. We follow the procedure in \cite{1}. Fix an interior point $\x_r$ of $O_{(\0, r)}$. Any ray starting at this point will intersect an obstacle other than $O_{(\0, r)}$ say $O_{(\k,t)}$ (This is due to Lemma 2.1 of \cite{1}).

Denote by $V_{(\k,t)}$ the set of directions in which any ray starting at $\x_r$ intersects the interior of $O_{(\k, t)}$. This set is open. Then the union of all the sets of this kind will be an open cover of the unit sphere which is compact. Hence, it has a finite sub-cover. This implies that there is a constant $M>0$ such that all rays of length $M$ starting at $\x_r$ intersect the interior of some $O_{(\k,t)}$, where $O_{(\k, t)}\neq O_{(\0, r)}$. \endofproof

Clearly, the set of vertices of $G$ in the case considered in this paper is much smaller than the set of vertices for the case of a torus with one obstacle.

We will denote  by $U$ {\it{the set of all ordered pairs of the form}} $(\u, r)$, where $\u$ is a unit vector in $\mathbb{Z}^m$, a vector with one component $\pm 1$ and the other components equal to 0, $r= 1,2,...,N$. Note that all elements of $U$ are vertices by condition 2.1.
\begin{lemma}
\label{8}
For any vertex $(\v, r)$ in the graph $G$, where $\v$ is in $\mathbb{Z}^m$, there is $(\u, r)$ in $U$ such that there is an edge from $(\v, r)$ to $(\u, r)$.
\end{lemma}
\emph{Proof}. Follows that in \cite{1}. Assume that $(\v, r)$ is a vertex, and $\v$ does not equal $\0$. Hence, there is a unit vector $\u\in\mathbb{Z}^m$ such that $\langle \u, \v\rangle\leq 0$. We show that we can separate $(\v, r)$ from $(\0, r)$ and $(\u+\v, r)$. 

Consider the triangle with vertices $A=\0$, $B=\v$ and $C=\u+\v$. Let $D$ and $E$ be the points on the sides $BA$, $BC$ respectively, whose distance from $B$ is $\frac{1}{2}$. Let $L$ be the straight line through $D$ and $E$. Since the angle at the vertex $B$ is at most $\frac{\pi}{2}$, then the distance of $B$ from $L$ is at least $\frac{\sqrt 2}{4}$. Since $\mid AD\mid\geq\mid BD\mid$ and $\mid CE\mid\geq\mid BE\mid$, the distances of $A$ and $C$ from $L$ are at least as large as the distance of $B$ from $L$. Now if we shift this triangle such that its vertices $A, \, B, \, C$ will be in the centers of the obstacles $(\0, r)$,  $(\v, r)$,  $(\u+\v, r)$ respectively, then the hyperplane of dimension $m-1$ through $L$ will separate $(\v, r)$ from $(\0, r)$ and $(\u+{\v}, r)$. i.e $(\v, r)$ is not between $(\0, r)$ and $(\u+\v, r)$. Thus, by the definition of an edge in $G$, there is an edge from $(\v, r)$ to $(\u, r)$. \endofproof

Note that, due to Corollary \ref{6}, there is also an edge from  $(\u, r)$ to  $(\v, r)$.

Now we can proof the main result of this section, the connectivity of $G$. 
\begin{theorem}
\label{9}
The graph $G$ is connected. In addition, for any two vertices $(\v, s)$ and $(\k, t)$ of $G$, there is a path of length at most $5$ from $(\v, s)$ to $(\k, t)$ in $G$ via elements of $U$.
\end{theorem}
\emph{Proof}. 
By Lemma \ref{8}, there are $(\u, s),\, (\n, t)$ in $U$ such that there are edges between them and $(\v, s)$, $(\k, t)$, respectively. From condition 2.1 and the definition of an edge in $G$, we can see that there are edges from $(\u, s)$ to $(\0, s)$, from $(\n, t)$ to $(\0, t)$, and from $(\0, s)$ to  $(\0, t)$.   Hence,  $(\v, s)$,  $(\u, s)$,  $(\0, s)$, $(\0, t)$, $(\n, t)$ $(\k, t)$ is a path of length $5$, and if  $(\u, s)= (\n, t)$ then  $(\v, s)$,  $(\u, s)$, $(\k, t)$ is a path of length $2$, by the note after Lemma \ref{8}. \endofproof

\section{The Admissible Rotation Set}
\setcounter{equation}{0}
In this section, we will prove our main results: the admissible rotation set is convex; and the rotation vectors of periodic trajectories of admissible type are dense in the admissible rotation set. Again we will use ideas from \cite{1} to prove these theorems. There are significant differences in the proof of the first theorem. For the second theorem, the proof is omitted as it is the same as the one in \cite{1} with minor change, namely, the length of the admissible sequence $C$ used in \cite{1} will be different in our case.

\begin{theorem}
The admissible rotation set $AR$ is convex.
\end{theorem}
\emph{Proof}. Fix $\u, \v$ in $AR$. We want to show that, for $t\in (0,1)$, $t\u+(1-t)\v$ is in $AR$. Fix $\epsilon>0$.

Since $\u$ and $\v$ are in $AR$, there are two admissible sequences $A$, $B$ and two trajectory pieces $T$, $S$ of type $A$ and $B$, respectively, such that 
\begin{equation} \norm[\bigg]{\frac{d(T)}{\mid T\mid}-\u}\leq \frac {\epsilon}{6s}\hspace{0.5 cm} \mbox{and}\hspace{0.5 cm}\norm[\bigg]{\frac{d(S)}{\mid S\mid}-\v}\leq\frac{\epsilon}{6(1-s)},\end{equation}
where $s=f(x)$, $$f(x)=\frac{\mid T\mid x}{\mid T\mid x + \mid S\mid (1-x)},$$ $x=\frac{p}{q}$ and $p,q$ are in $\mathbb{Z}^+$. Clearly, this function of $x$ is continuous on $[0,1]$, takes values $0$ at $0$ and $1$ at $1$. Hence, the image of the set of rational numbers is dense in $[0,1]$. Therefore, we can approximate $t$ by $s$ with an arbitrary accuracy. Thus, we can choose values for $p$ and $q$ such that \begin{equation}\norm{t-s}\leq\frac{\epsilon}{6k},\, k=max\{\norm{\u},\norm{\v}\}.\end{equation}
Now, by Lemma \ref{8}, there are admissible sequences $C_1, C_2, C_3$ of length at most $5$ via elements of $U$ such that 
$$D=AC_1AC_1......AC_1AC_2BC_3BC_3.....BC_3B$$
is admissible, and so there is an admissible trajectory piece $Q$ of type $D$ by Theorem \ref{1}. Assume that the block $A$ appears $p$ times on $D$ and $B$ appears $q-p$ times, where $q>p$. Let $d_A$ represents the displacement due to the repetition of the block $A$, $d_B$ represents the displacement due to the repetition of the block $B$. Then, as in \cite{1},
$$\norm{d_A-pd(T)}\leq 2pd\hspace{0.5 cm} \mbox{and} \hspace{0.5 cm}\norm{d_B-(q-p)d(S)}\leq2(q-p)d,$$ where $d$ is maximal diameter of the obstacles.
In addition, the displacement due to each $C_i$ is at most $5+2d$, so the total displacement due to all those blocks is at most $(q-1)(5+2d)<q(5+2d)$.

Thus \begin{equation}\norm{d(Q)-pd(T)-(q-p)d(S)}\leq2pd+2(p-q)d+q(5+2d)=4qd+5q.\end{equation}
Similarly, by Note \ref{4}, \begin{equation} \mid\mid Q\mid-p\mid T\mid-(q-p)\mid S\mid\mid\leq4qd+5q.\end{equation}
Now 
\begin{eqnarray}
&&\norm[\bigg]{\frac{d(Q)}{\mid Q\mid}-s\frac{d(T)}{\mid T\mid}-(1-s)\frac{d(S)}{\mid S\mid}}\nonumber\\
&=&\norm[\bigg]{\frac{d(Q)}{\mid Q\mid}-\frac{pd(T)}{p\mid T\mid+(q-p)\mid S\mid}-\frac{(q-p)d(S)}{p\mid T\mid+(q-p)\mid S\mid}}\nonumber\\
&\leq&\norm[\bigg]{\frac{d(Q)(p\mid T\mid +(q-p)\mid S\mid)-d(Q)\mid Q\mid}{\mid Q\mid(p\mid T\mid +(q-p)\mid S\mid)}}\nonumber\\
&&+\norm[\bigg]{\frac{d(Q) \mid Q\mid-pd(T)\mid Q\mid-(q-p)d(S)\mid Q\mid}{\mid Q\mid(p\mid T\mid+(q-p)\mid S\mid)}}\nonumber\\
&\leq&\frac{d(Q)}{\mid Q\mid}\frac{(4qd+5q)}{p\mid T\mid+\mid S\mid (q-p)}+\frac{4qd+5q}{p\mid T\mid+\mid S\mid(q-p)}\nonumber\\
&\leq&\frac{8qd+10q}{p\mid T\mid+\mid S\mid (q-p)}=\frac{8d+10}{\frac{p}{q}\mid T\mid+\mid S\mid(1-\frac{p}{q})}\hspace{3cm}\mbox{(as $\frac{d(Q)}{\mid Q\mid}\leq 1$)}.\nonumber
\end{eqnarray}
For any $p$ and $q$ either $\frac{p}{q}$ or $(1-\frac{p}{q})\geq\frac{1}{2}$. If $\frac{p}{q}\geq\frac{1}{2}$ then 
\begin{equation}\frac{8d+10}{\frac{p}{q}\mid T\mid+\mid S\mid(1-\frac{p}{q})}\leq\frac{8d+10}{\frac{1}{2}\mid T\mid}=\frac{16d+20}{\mid T\mid}\leq\frac{\epsilon}{3}\end{equation}
when we take $\mid T\mid$ large enough. We get a similar result when $(1-\frac{p}{q}) \geq\frac{1}{2}$.

Thus, 
\begin{eqnarray}
&&\norm[\bigg]{\frac{d(Q)}{\mid Q\mid}-t\u-(1-t)\v}\nonumber\\
&\leq&\norm[\bigg]{\frac{d(Q)}{\mid Q\mid}-s\u-(1-s)\v}+\norm[\bigg]{s\u+(1-s)\v-t\u-(1-t)\v}\nonumber\\
&\leq&\norm[\bigg]{\frac{d(Q)}{\mid Q\mid}-s\frac{d(T)}{\mid T\mid}-(1-s)\frac{d(S)}{\mid S\mid}}+\norm[\bigg]{s\frac{d(T)}{\mid T\mid}-s\u}\nonumber\\
&&+\norm[\bigg]{(1-s)\frac{d(S)}{\mid S\mid}-(1-s)\v}+2\norm{s-t}k\hspace{4cm}\mbox{(by $( 4.2)$)}\nonumber\\
&\leq&\epsilon. \hspace{8 cm}\mbox{(by $(4.1), (4.2)$ and $(4.5)$)}\nonumber
\end{eqnarray}
This completes the proof. \endofproof
\begin{theorem}
The rotation vectors of periodic orbits of admissible type are dense in the admissible rotation set.
\end{theorem}
\emph{Proof}. The proof in \cite{1} is applied here. We just need to change the length of the the admissible sequence $C$ used in \cite{1} to $5$.  \endofproof

The following theorem is similar to Theorem 4.15 in \cite{1}.
\begin{theorem}
The admissible rotation set is a proper subset of the general rotation set.
\end{theorem}

\emph{Proof}. The proof in \cite{1} shows that the admissible rotation set is contained in the unit open ball. Here, we provide an example which shows that the vector (1,0,0,0,....,0) is in the general rotation set. Consider the trajectory piece joining the lowest point $\x_k$ in the boundary of the  lowest obstacle $O_{(\0,r_0)}$ in the hyper cube with the leading vector $\0$ with the highest point $y_k$ in the boundary of the highest obstacle in the hypercube with leading vector $(k, -1, 0,0,...,0)$, where $k$ is large enough so that the trajectory piece does not intersect any obstacles other than those. Then as $k$ goes to infinity the ratio $\frac{y_k-x_k}{t_k}$ goes to $(1,0,0,...,0)$. 

\footnotesize
{\bf Acknowledgements.} Financial support by the Ministry of Higher Education in Saudi Arabia is gratefully acknowledged.

\end{document}